\documentclass[12pt,a4paper]{article}
\usepackage{amsmath,amsthm,amssymb}
\begin{document}

%\input{setsymb}

%%
%  Symbols for set theory
%
%
%  Ideals on the real line
%
\def\N{\mathcal{N}}
\def\M{\mathcal{M}}
%
%  Set-theoretic operations/relations
%
\def\domn{\leq^*}
\def\ssm{\smallsetminus}
\def\restrictedto{\!\upharpoonright\!}
\def\op#1#2{\langle #1,#2\rangle}
\def\conc{{}^\frown\!}			% defined by Roslanowski
\def\st{:}
\def\card#1{\lvert#1\rvert}
\def\size#1{\lvert#1\rvert}
\def\lh#1{\lvert#1\rvert}
\def\Seq#1{\left\langle{#1}\right\rangle}
%
%  Logical symbols
%
\def\forallbutfin{\forall^{\infty}}
\def\existinf{\exists^{\infty}}
%
%  Forcing
%
\newcommand{\forces}[1][{}]{\Vdash_{#1}}
\newcommand{\forcestext}[2][{}]{\Vdash_{#1}\text{``}{#2}\text{''}}
\def\supp{\operatorname{\textup{\textsf{supp}}}}
\newcommand\tvalue[2][{}]{{[\![#2]\!]_{#1}}}
\def\V{\mathbf{V}}
\def\W{\mathbf{W}}
\def\incrfs{\omega^{\uparrow\omega}}
\def\incrseqs{\omega^{\uparrow<\omega}}

\def\pdomn{\sqsubseteq}
\def\rk#1{\operatorname{rank}(#1)}

%\input{thmdef}

%
%  Theorem style definitions
%
\theoremstyle{plain}
\newtheorem{thm}{Theorem}[section]
\newtheorem{lem}[thm]{Lemma}
\newtheorem{cor}[thm]{Corollary}
\newtheorem{prop}[thm]{Proposition}
\newtheorem{fact}[thm]{Fact}
\newtheorem{claim}{Claim}
\newtheorem{conj}[thm]{Conjecture}
\newtheorem{q}[thm]{Question}
\theoremstyle{definition}
\newtheorem{defn}[thm]{Definition}
\theoremstyle{remark}
\newtheorem{ack}{Acknowledgement}
	\renewcommand{\theack}{}
\newtheorem{remark}{Remark}

\def\bP{\mathbb{P}}
\def\bC{\mathbb{C}}

\title{Hechler's theorem for the meager ideal}
\author{Tomek Bartoszy\'nski\thanks{%
	The first author was  partially supported by 
	NSF grant DMS 0200671.}
	\ and Masaru Kada\thanks{%
	The second author was supported by 
	Grant-in-Aid for Young Scientists (B) 14740058, JSPS.
	}
}
%\date{August 8, 2002}
\date{}
\maketitle

%%%%%%%%
\renewcommand{\thefootnote}{}
\renewcommand{\thefootnote}{\fnsymbol{footnote}}
%%%%%%%%

\begin{abstract}
We prove the following theorem: 
For a partially ordered set $Q$ such that every countable subset has 
a strict upper bound, 
there is a forcing notion satisfying ccc such that, 
in the forcing model, 
there is a basis of the meager ideal of the real line 
which is order-isomorphic to $Q$ with respect to set-inclusion. 
This is a variation of Hechler's classical result in the theory of forcing. 
\end{abstract}

\section{Introduction}

For $f,g\in\omega^\omega$, 
we say $f\domn g$ if $f(n)\leq g(n)$ for all but finitely many $n<\omega$. 
The following theorem, which is due to Hechler \cite{He:cof}, 
is a classical result in the theory of forcing 
(See also \cite{Bur:hechler}). 

\begin{thm}
Suppose that $(Q,\leq)$ is a partially ordered set such that 
every countable subset of $Q$ has a strict upper bound in $Q$, 
that is, 
for any countable set $A\subseteq Q$ there is $b\in Q$ 
such that $a<b$ for all $a\in A$. 
Then there is a forcing notion $\bP$ satisfying ccc 
such that, 
in the forcing model by $\bP$, 
$(\omega^\omega,\domn)$ contains a cofinal subset 
$\{f_a\st a\in Q\}$ which is order-isomorphic to $Q$, 
that is, 
\begin{enumerate}
\item for every $g\in\omega^\omega$ 
	there is $a\in Q$ such that $g\domn f_a$, and 
\item for $a,b\in Q$, $f_a\domn f_b$ if and only if $a\leq b$. 
\end{enumerate}
\end{thm}

Soukup \cite{So:specnote} 
asked if the statement of Hechler's theorem holds for the meager ideal 
or the null ideal of the real line with respect to set-inclusion. 

In this paper 
we give a positive answer for the meager ideal. 
The basic idea of the construction of the forcing notion 
is the same as Hechler's original proof, 
but we modify it to fit in our context. 

The question for the null ideal was answered positively 
by the second author~\cite{Kada:embedn}.

Let $\incrfs$ and $\incrseqs$ be 
the set of strictly increasing functions in $\omega^\omega$ 
and the set of strictly increasing sequences in $\omega^{<\omega}$ 
respectively. 
For $f,g\in\incrfs$, 
$f\pdomn g$ if for all but finitely many $n<\omega$ 
there is $k<\omega$ such that 
$[f(k),f(k+1))\subseteq[g(n),g(n+1))$. 
We say 
$d\in\incrfs$ is a $\pdomn$-dominating real over a model $\V$ of ZFC 
if $f\pdomn d$ for all $f\in\incrfs\cap\V$. 

For $x\in 2^\omega$ and $f\in\incrfs$, 
define a meager set $E_{x,f}\subseteq 2^\omega$ by the following: 
\[
E_{x,f}=\{z\in 2^\omega\st
	\exists m<\omega\,\forall n\geq m\,\exists j\in[f(n),f(n+1))\,
	(z(j)\neq x(j))\}.
\]

\begin{lem}\label{lem:pdomn}
For $x\in 2^\omega$ and $f,g\in\incrfs$, 
if $f\pdomn g$ then $E_{x,f}\subseteq E_{x,g}$. 
\end{lem}

\begin{proof}
Clear. 
\end{proof}

\begin{lem}\label{lem:notpdomn}
For $x,y\in 2^\omega$ and $f,g\in\incrfs$, 
if $f\not\pdomn g$ then $E_{x,f}\not\subseteq E_{y,g}$. 
\end{lem}

\begin{proof}
Suppose that $x,y\in 2^\omega$, $f,g\in\incrfs$ and $f\not\pdomn g$. 
Let 
$A=\{n<\omega\st
	[f(k),f(k+1))\not\subseteq[g(n),g(n+1))\text{ for all }k<\omega\}$. 
By the assumption, $A$ is an infinite subset of $\omega$. 
Define $z\in 2^\omega$ as follows: 
\[
z(j)=
\begin{cases}
y(j)	&	j\in[g(n),g(n+1))\text{ for some }n\in A \\
1-x(j)	&	\text{otherwise}
\end{cases}.
\]
It is easy to see that $z\in E_{x,f}\ssm E_{y,g}$. 
\end{proof}

\begin{lem}\label{lem:cover}
Suppose that $\V$ is a model of \textup{ZFC}, 
$c$ is a Cohen real over $\V$, 
and $d$ is a $\pdomn$-dominating real over $\V[c]$. 
Then, for any Borel meager set $X\subseteq 2^\omega$ 
which is coded in $\V$, 
we have $X\subseteq E_{c,d}$. 
\end{lem}

\begin{proof}
Fix $x\in X$. 
Since $X$ is coded in $\V$ 
and $c$ is a Cohen real over $\V$, 
there are infinitely many $j<\omega$ 
such that $x(j)\neq c(j)$. 
We can define an infinite set 
$D_x=\{j<\omega\st x(j)\neq c(j)\}\subseteq\omega$ 
in $\V[c]$. 
Since $d$ is a $\pdomn$-dominating real over $\V[c]$, 
for all but finitely many $n<\omega$ we have 
$D_x\cap[d(n),d(n+1))\neq\emptyset$, 
and hence $x\in E_{c,d}$. 
\end{proof}

We will use the following standard fact about partially ordered sets. 
See \cite{Sz:poset} for the proof. 

\begin{prop}\label{prop:poset}
If $(P,\leq)$ is a partially ordered set and $c\in P$, 
then the partial order $\leq$ on $P$ 
can be extended to a linear order $\leq'$ so that 
$c\leq' y$ for every $y\in P$ which is $\leq$-incomparable to $c$. 
\end{prop}

\begin{remark}
Note that 
$f\pdomn g$ is 
neither a sufficient nor a necessary condition for $f\domn g$.
We say 
$d\in\omega^\omega$ is a $\domn$-dominating real over $\V$ 
if $f\domn d$ for all $f\in\omega^\omega\cap\V$. 
It is easy to see that 
a $\pdomn$-dominating real over $\V$ is also $\domn$-dominating over $\V$, 
but the converse does not hold in general. 
However, we can construct a $\pdomn$-dominating real 
from a $\domn$-dominating real 
(See \cite[Theorem~2.10]{Bl:combchar} for the proof). 
\end{remark}

\section{The main theorem}

Let $(Q,\leq)$ be a partially ordered set such that 
every countable subset of $Q$ has a strict upper bound in $Q$, 
that is, 
for any countable set $A\subseteq Q$ there is $b\in Q$ 
such that $a<b$ for all $a\in A$. 
Extend the order to $Q^*=Q\cup\{Q\}$ 
by letting $a<Q$ for each $a\in Q$. 
Let $R\subseteq Q$ be a well-founded cofinal subset. 
Define the rank function on the well-founded set $R^*=R\cup\{Q\}$ 
in the usual way. 
For $a\in Q\ssm R$, 
let $\rk{a}=\min\{\rk{b}\st b\in R^*\text{ and }a<b\}$. 
For $x,y\in Q^*$, 
we say $x\ll y$ if $x<y$ and $\rk{x}<\rk{y}$. 
For $x\in Q^*$, 
let $Q_x=\{y\in Q\st y\ll x\}$. 
For $F\subseteq Q^*$, let $\bar{F}=\{\rk{x}\st x\in F\}$. 

Let $\bC=2^{<\omega}$ 
be the forcing notion adding one Cohen real. 

We define forcing notions $\bP_a$ by induction on $\rk{a}$ for $a\in Q^*$. 

%For $a\in Q^*$, 
A condition of a forcing notion $\bP_a$ 
is of the form 
$p=(\{s_\alpha\st\alpha\in\bar{F}\},\{(t_b,\dot{f}_b)\st b\in F\})$ 
with the following: 
\begin{enumerate}
\item $F$ is a finite subset of $Q_a$; 
\item For $\alpha\in\bar{F}$, $s_\alpha\in\bC$; 
\item For $b\in F$, 
	$t_b\in\incrseqs$, and 
	$\dot{f}_b$ is a $\bP_b*\bC$-name 
	for a function in $\incrfs$. 
\end{enumerate}

For $p\in\bP_a$ and $b<a$, 
define $p\restrictedto b\in\bP_b$ 
by letting 
$p\restrictedto b
	=(\{s_\alpha\st\alpha\in\bar{F_b}\},\{(t_c,\dot{f}_c)\st c\in F_b\})$ 
where $F_b=F\cap Q_b$. 

For conditions 
$p=(\{s^p_\alpha\st\alpha\in\bar{F^p}\},\{(t^p_b,\dot{f}^p_b)\st b\in F^p\})$ 
and 
$q=(\{s^q_\alpha\st\alpha\in\bar{F^q}\},\{(t^q_b,\dot{f}^q_b)\st b\in F^q\})$ 
in $\bP_a$, 
$p\leq q$ if the following hold: 

\begin{enumerate} 
\item $F^q\subseteq F^p$; 
\item For $\alpha\in\bar{F^q}$, $s^q_\alpha\subseteq s^p_\alpha$;
\item For $b\in F^q$, 
	$t^q_b\subseteq t^p_b$, and the condition 
	$\op{p\restrictedto b}{s^p_\beta}\in\bP_b*\bC$ 
	forces that: 
	\begin{enumerate}
	\item for all $n<\omega$ there is $k<\omega$ 
		such that $[\dot{f}^q_b(k),\dot{f}^q_b(k+1))
			\subseteq[\dot{f}^p_b(n),\dot{f}^p_b(n+1))$, and 
	\item for all $n\in\lh{t^p_b}\ssm\lh{t^q_b}$ 
		there is $k<\omega$ such that 
		$[\dot{f}^q_b(k),\dot{f}^q_b(k+1))
			\subseteq[t^p_b(n-1),t^p_b(n))$, 
	\end{enumerate}
where $\beta=\rk{b}$; 
\item For $b,c\in F^q$, 
	if $b<c$ and $\rk{b}=\rk{c}$, 
	then for all $n\in\lh{t^p_c}\ssm\lh{t^q_c}$ 
	there is $k<\lh{t^p_b}$ 
	such that $[t^p_b(k-1),t^p_b(k))\subseteq[t^p_c(n-1),t^p_c(n))$. 
\end{enumerate}

\begin{lem}
$\bP_Q$ satisfies ccc. 
\end{lem}

\begin{proof}
A standard $\Delta$-system argument. 
\end{proof}

\begin{lem}
For $a,b\in Q^*$ with $a\ll b$, 
the inclusion from $\bP_a$ to $\bP_b$ is a complete embedding. 
\end{lem}

\begin{proof}
Clear. 
\end{proof}

Let $\V$ be the ground model, and $G$ be a $\bP_Q$-generic filter over $\V$. 
For $a\in Q$, let $G\restrictedto a=\{p\restrictedto a\st p\in G\}$. 

Work in $\V[G]$. 
We assume that 
each $p\in\bP_Q$ 
is represented as 
$p=(\{s^p_\alpha\st\alpha\in\bar{F^p}\},\{(t_a^p,\dot{f}_a^p)\st a\in F^p\})$. 
%
%$p^c$ and $p^d$ denote the 
%first and second coordinate of $p$ respectively. 
%
%we denote $F^p$, $s_\alpha^p$, $t_a^p$ and $\dot{f}_a^p$ 
%for respective components of $p$. 
%
For $\alpha<\rk{Q}$, 
let $c_\alpha=\bigcup\{s_\alpha^p\st
		p\in G\text{ and }\alpha\in\bar{F^p}\}$, 
and 
for $a\in Q$, 
let $d_a=\bigcup\{t_a^p\st
		p\in G\text{ and }a\in F^p\}$. 

Clearly, if $a\in Q$ and $\alpha=\rk{a}$, 
then $c_\alpha\in 2^\omega$ and 
$c_\alpha$ is a Cohen real over $\V[G\restrictedto a]$. 

\begin{lem}\label{lem:total}
For each $a\in Q$, 
$d_a\in\incrfs$, 
that is, $d_a$ is defined on all of $\omega$.  
\end{lem}

\begin{proof}
We will show that, 
for every 
$q=(\{s^q_\alpha\st\alpha\in\bar{F^q}\},\{(t^q_b,\dot{f}^q_b)\st b\in F^q\})
	\in\bP_Q$ 
there is 
$p=(\{s^p_\alpha\st\alpha\in\bar{F^p}\},\{(t^p_b,\dot{f}^p_b)\st b\in F^p\})
	\leq q$ 
such that for every $b\in F^q$ we have $\lh{t^p_b}>\lh{t^q_b}$. 

Let $q\in\bP_Q$ and $\alpha=\max\bar{F^q}$. 
We work by induction on $\alpha$. 

Let 
$q_{<\alpha}=(\{s^q_\beta\st\beta\in\bar{F^q}\cap\alpha\},
	\{(t^q_b,\dot{f}^q_b)\st b\in F^q\text{ and }\rk{b}<\alpha\})$. 
It is easily seen that $q_{<\alpha}\in\bP_Q$. 
By the induction hypothesis, 
there is a condition 
$r=(\{s^r_\beta\st\beta\in\bar{F^r}\},\{(t^r_b,\dot{f}^r_b)\st b\in F^r\})
	\leq q_{<\alpha}$ 
such that for every $b\in F^q$ 
if $\rk{b}<\alpha$ then $\lh{t^r_b}>\lh{t^q_b}$. 

Define 
$p^0=
(\{s^{p^0}_\beta\st\beta\in\bar{F^{p^0}}\},
	\{(t^{p^0}_b,\dot{f}^{p^0}_b)\st b\in F^{p^0}\})$ 
as follows: 
\begin{enumerate}
\item $F^{p^0}=\{b\in F^r\st\rk{b}<\alpha\}\cup F^q$; 
\item For $\beta\in\bar{F^{p^0}}\cap\alpha$, $s^{p^0}_\beta=s^r_\beta$; 
\item For $b\in F^{p^0}$ with $\rk{b}<\alpha$, 
	$t^{p^0}_b=t^r_b$ and $\dot{f}^{p^0}_b=\dot{f}^r_b$; 
\item $s^{p^0}_\alpha=s^q_\alpha$;
\item For $b\in F^{p^0}$ with $\rk{b}=\alpha$, 
	$t^{p^0}_b=t^q_b$ and $\dot{f}^{p^0}_b=\dot{f}^q_b$. 
\end{enumerate}

It is easy to check 
that $p^0\in\bP_Q$ and $p^0\leq q$. 

Extend the order $<$ on $\{b\in F^{p^0}\st\rk{b}=\alpha\}$ 
to a linear order $<'$, 
say $\{b\in F^{p^0}\st\rk{b}=\alpha\}=\{b_1,\ldots,b_n\}$ 
with $b_1<'\cdots<'b_n$. 
We will inductively define conditions $p^i$ for $i=1,\ldots,n$ 
such that $p^0\geq p^1\geq\cdots\geq p^n$. 

Suppose that $1\leq i\leq n$ and 
$p^j=
(\{s^{p^j}_\beta\st\beta\in\bar{F^{p^j}}\},
	\{(t^{p^j}_b,\dot{f}^{p^j}_b)\st b\in F^{p^j}\})$ 
is already defined for $j\leq i-1$. 
To find $p^i\leq p^{i-1}$, we will construct a decreasing sequence 
$p^{i-1}=r^0\geq r^1\geq\cdots\geq r^{2^i-1}=p^i$, 
where $r^k=(\{s^{r^k}_\beta\st\beta\in\bar{F^{r^k}}\},
	\{(t^{r^k}_b,\dot{f}^{r^k}_b)\st b\in F^{r^k}\})$ 
for $0\leq k\leq 2^i-1$, 
in the following way. 

\textit{Step 1.}
Find 
$w=(\{s^w_\beta\st\beta\in\bar{F^w}\},\{(t^w_b,\dot{f}^w_b)\st b\in F^w\})
	\in\bP_{b_1}$
with 
$w\leq r^0\restrictedto b_1$, 
$v\leq s^{r^0}_\alpha$ 
and $h\in\incrseqs$ 
such that $h(\lh{h}-2)\geq t^{r^0}_{b_1}(\lh{t^{r^0}_{b_1}}-1)$ 
and 
$(w,v)\forcestext[\bP_{b_1}*\bC]{h\subseteq\dot{f}^{r^0}_{b_1}}$. 
Define $r^1$ as follows: 
\begin{enumerate}
\item $F^{r^1}=F^w\cup F^{r^0}$; 
\item For $\beta\in\bar{F^w}$, $s^{r^1}_\beta=s^w_\beta$; 
\item For $b\in F^w$, 
	$t^{r^1}_b=t^w_b$ 
	and $\dot{f}^{r^1}_b=\dot{f}^w_b$;
\item $s^{r^1}_\alpha=v$;
\item $t^{r^1}_{b_1}=t^{r^0}_{b_1}\conc\Seq{l}$, 
	where $l\geq\max(\{t^{r^0}_{b_i}(\lh{t^{r^0}_{b_i}}-1)
		\st 1\leq i\leq n\}
		\cup\{h(\lh{h}-1)\})$, 
	and $\dot{f}^{r^1}_{b_1}=\dot{f}^{p^{r^0}}_{b_1}$; 
\item For $\beta\in\bar{F^{r^1}}\ssm(\bar{F^w}\cup\{\alpha\})$, 
	$s^{r^1}_\beta=s^{r^0}_\beta$; 
\item For $b\in F^{r^1}\ssm (F^w\cup\{b_1\})$, 
	$t^{r^1}_b=t^{r^0}_b$ 
	and $\dot{f}^{r^1}_b=\dot{f}^{p^{r^0}}_b$. 
\end{enumerate}

It is easy to check that $r^1\in\bP_Q$ and $r^1\leq r^0$. 

\textit{Step 2.}
Again we find $w\leq r^1\restrictedto b_1$, $v\leq s^{r^1}_\alpha$ and $h$ 
such that 
$h(\lh{h}-2)\geq t^{r^1}_{b_1}(\lh{t^{r^1}_{b_1}}-1)$ 
and 
$(w,v)\forcestext[\bP_{b_1}*\bC]{h\subseteq\dot{f}^{r^1}_{b_1}}$, 
and form a condition $r^2\leq r^1$ as in the previous step. 

\textit{Step 3.}
Now we look at $b_2$. 
Find $w\leq r^1\restrictedto b_2$, $v\leq s^{r^2}_\alpha$ and $h$ 
such that 
$h(\lh{h}-2)\geq t^{r^2}_{b_2}(\lh{t^{r^2}_{b_2}}-1)$ 
and 
$(w,v)\forcestext[\bP_{b_2}*\bC]{h\subseteq\dot{f}^{r^2}_{b_2}}$. 
Then we extend $t^{r^2}_{b_2}$ 
and form a condition $r^3$ in the same way. 

\setlength{\unitlength}{1mm}
\begin{figure}
\begin{picture}(130,65)(-15,-5)

\put(-10,50){\makebox(0,0){$b_1$}}
\put(-10,40){\makebox(0,0){$b_2$}}
\put(-10,30){\makebox(0,0){$b_3$}}
\put(-10,20){\makebox(0,0){$\vdots$}}
\put(-10,10){\makebox(0,0){$b_{n-1}$}}
\put(-10,0){\makebox(0,0){$b_n$}}

\thinlines
\put(0,50){\vector(1,0){110}}
\put(0,40){\vector(1,0){110}}
\put(0,30){\vector(1,0){110}}
\put(0,10){\vector(1,0){110}}
\put(0,0){\vector(1,0){110}}

\put(55,20){\makebox(0,0){$\vdots$}}

\put(15,50){\circle*{1.5}}
\put(10,40){\circle*{1.5}}
\put(17,30){\circle*{1.5}}
\put(13,10){\circle*{1.5}}
\put(30,0){\circle*{1.5}}

\thicklines
\put(0,50){\line(1,0){15}}
\put(0,40){\line(1,0){10}}
\put(0,30){\line(1,0){17}}
\put(0,10){\line(1,0){13}}
\put(0,0){\line(1,0){30}}

\put(15,51){\makebox(0,0)[b]{$t^{r^0}_{b_1}(\lh{t^{r^0}_{b_1}}-1)$}}
\put(10,41){\makebox(0,0)[b]{$t^{r^0}_{b_2}(\lh{t^{r^0}_{b_2}}-1)$}}
\put(17,31){\makebox(0,0)[b]{$t^{r^0}_{b_3}(\lh{t^{r^0}_{b_3}}-1)$}}
\put(13,11){\makebox(0,0)[b]{$t^{r^0}_{b_{n-1}}(\lh{t^{r^0}_{b_{n-1}}}-1)$}}
\put(30,1){\makebox(0,0)[b]{$t^{r^0}_{b_n}(\lh{t^{r^0}_{b_n}}-1)$}}

\put(36,50){\circle*{1}}
\put(40,50){\circle*{1}}
\put(44,40){\circle*{1}}
\put(48,50){\circle*{1}}
\put(52,50){\circle*{1}}
\put(56,40){\circle*{1}}
\put(60,30){\circle*{1}}
\put(64,50){\circle*{1}}
\put(80,10){\circle*{1}}
\put(96,10){\circle*{1}}
\put(100,0){\circle*{1}}

\put(36,51){\makebox(0,0)[b]{$1$}}
\put(40,51){\makebox(0,0)[b]{$2$}}
\put(44,41){\makebox(0,0)[b]{$3$}}
\put(48,51){\makebox(0,0)[b]{$4$}}
\put(52,51){\makebox(0,0)[b]{$5$}}
\put(56,41){\makebox(0,0)[b]{$6$}}
\put(60,31){\makebox(0,0)[b]{$7$}}
\put(64,51){\makebox(0,0)[b]{$8$}}
\put(72,51){\makebox(0,0)[b]{$\cdots$}}
\put(74,41){\makebox(0,0)[b]{$\cdots$}}
\put(76,31){\makebox(0,0)[b]{$\cdots$}}
\put(80,11){\makebox(0,0)[b]{$2^n{-}1$}}
\put(96,11){\makebox(0,0)[b]{$2^{n+1}{-}2$}}
\put(100,1){\makebox(0,0)[b]{$2^{n+1}{-}1$}}
\end{picture}
\caption{Extending $t$'s at $b_1,\ldots,b_n$}
\label{fig:extension}
\end{figure}
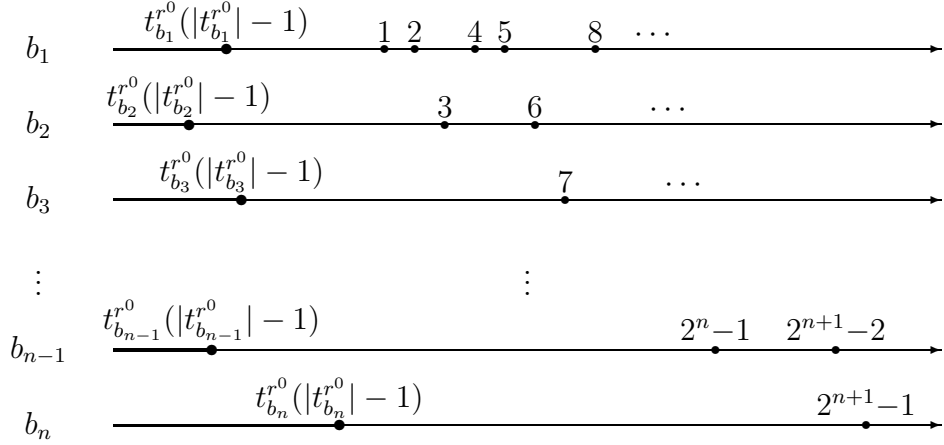

We extend $t$'s for $b_i$'s and define conditions $r^k$'s in the same way, 
step by step 
along the order shown in Figure~\ref{fig:extension}. 
That is, 
for each step, we set the value of the first open place of $t$ at $b_i$ 
so that it exceeds all values which are already set, 
and the last interval determined by $t$ contains some interval 
determined by the corresponding $\dot{f}$.  
If two consecutive values of $t$ at $b_i$ are defined, 
then we extend $t$ at $b_{i+1}$. 
In the $(2^{n+1}-1)$-st step, we can extend $t$ at $b_{n}$ 
and then every sequence has been extended. 

Finally, let $p=p^n$. 
It is straightforward to check that $p$ is as desired. 
\end{proof}

Now it is easy to see that, 
for $a\in Q$ with $\rk{a}=\alpha$, 
$d_a$ is a $\pdomn$-dominating real over $\V[G\restrictedto a][c_\alpha]$. 

%Throughout this paper, 
%we use the notation $c_\alpha$ and $d_a$ 
%also to denote the corresponding $\bP_Q$-names. 

\begin{lem}\label{lem:incomp}
For $a,b\in Q$, 
$d_a\pdomn d_b$ if and only if $a\leq b$. 
\end{lem}

\begin{proof}
%Fix $a,b\in Q$. 
It is easy to see that $a\leq b$ implies $d_a\pdomn d_b$. 
Now assume $a\not\leq b$. 
Let $\alpha=\rk{a}$ and $\beta=\rk{b}$. 
We will show $d_a\not\pdomn d_b$ by a similar argument 
as in the proof of Lemma~\ref{lem:total}. 

Fix $N<\omega$ and 
$q=(\{s^q_\gamma\st\gamma\in\bar{F^q}\},\{(t^q_c,\dot{f}^q_c)\st b\in F^q\})
	\in\bP_Q$. 
We may assume that $a,b\in F^q$. 
Let $M=\max\{N,\lh{t^q_b}\}$. 
We will find a condition $p\leq q$ 
which forces that the interval $[d_b(M),d_b(M+1))$ does not contain 
any interval of the form $[d_a(k),d_a(k+1))$. 

Extend the order $<$ on 
$\{x\in F^q\st x\leq b\text{ and }\rk{x}=\beta\}$ 
to a linear order $<'$, 
say 
$\{x\in F^q\st x\leq b\text{ and }\rk{x}=\beta\}
	=\{x_1,\ldots,x_n\}$ 
with $x_1<'\cdots<'x_{n-1}<'x_n=b$. 
Using the method in the proof of Lemma~\ref{lem:total} 
along the order $<'$, 
extend $q$ to $q^*$ 
so that $\lh{t^{q^*}_b}=M+2$, $t^{q^*}_a=t^q_a$ and 
$t^{q^*}_b(M)>t^{q^*}_a(\lh{t^{q^*}_a}-1)$. 
This is possible because $a$ is not below $b$ 
and so $t^q_a$ is never extended through this process. 

Next, extend the order $<$ on 
$\{y\in F^{q^*}\st y\leq a\text{ and }\rk{y}=\alpha\}$ 
to a linear order $<'$, 
say $\{y\in F^{q^*}\st y\leq a\text{ and }\rk{y}=\alpha\}
	=\{y_1,\ldots,y_m\}$ 
with $y_1<'\cdots<'y_{m-1}<'y_m=a$. 
Again, as in the proof of Lemma~\ref{lem:total}, 
we extend $q^*$ to $p$ 
so that 
$\lh{t^p_a}=\lh{t^{q^*}_a}+1$ and 
$t^p_a(\lh{t^{q^*}_a})>t^p_b(M+1)=t^{q^*}_b(M+1)$. 

It is easy to check that $p$ forces that 
the interval $[d_b(M),d_b(M+1))$ contains no value of $d_a$, 
which concludes the proof. 
\end{proof}

For $a\in Q$, let $E_a=E_{c_\alpha,d_a}$ where $\alpha=\rk{a}$. 

\begin{lem}\label{lem:cohendomcover}
Let $a\in R$. 
If $X\subseteq 2^\omega$ is a Borel meager set 
which is coded in $\V[G\restrictedto a]$, 
then $X\subseteq E_a$. 
\end{lem}

\begin{proof}
Follows from Lemmata~\ref{lem:cover} and \ref{lem:total}. 
\end{proof}

\begin{cor}\label{cor:isomorphic}
For $a,b\in Q$, $a\leq b$ if and only if $E_a\subseteq E_b$. 
\end{cor}

\begin{proof}
If $a\ll b$, then $E_a\subseteq E_b$ by Lemma~\ref{lem:cohendomcover}. 
If $a\leq b$ and $\rk{a}=\rk{b}$, 
then $E_a\subseteq E_b$ by Lemmata~\ref{lem:pdomn} and \ref{lem:incomp}. 
If $a\not\leq b$, 
then $E_a\not\subseteq E_b$ follows from 
Lemmata~\ref{lem:notpdomn} and \ref{lem:incomp}. 
\end{proof}

\begin{cor}\label{cor:cofinal}
In $\V[G]$, for every meager set $X$ there is $a\in Q$ 
such that $X\subseteq E_a$.  
\end{cor}

\begin{proof}
By the definition of $\bP_Q$ and the assumption on $Q$, 
every Borel set in $\V[G]$ is coded in $\V[G\restrictedto a]$ 
for some $a\in Q$. 
\end{proof}

Now we have the following main theorem. 

\begin{thm}
Let $\M$ be the collection of meager sets in $2^\omega$. 
Suppose that $Q$ is a partially ordered set such that 
every countable subset of $Q$ has a strict upper bound in $Q$. 
Then in the forcing model by $\bP_Q$, 
$(\M,\subseteq)$ contains a cofinal subset 
$\{E_a\st a\in Q\}$ which is order-isomorphic to $Q$, 
that is, 
\begin{enumerate}
\item for every $X\in\M$ there is $a\in Q$ 
	such that $X\subseteq E_a$, and 
\item for $a,b\in Q$, $E_a\subseteq E_b$ if and only if $a\leq b$. 
\end{enumerate}
\end{thm}

\begin{remark}
The forcing $\bP_Q$ adds Cohen reals indexed by the ranks of $Q$ 
and dominating reals indexed by $Q$ itself. 
%If $Q$ is well-founded, 
%we can prove the main theorem by adding 
%both Cohen and dominating reals indexed by $Q$. 
One might add 
both Cohen and dominating reals indexed by $Q$, 
say $\{(c_a,d_a)\st a\in Q\}$, 
and set $E_a=E_{c_a,d_a}$ for $a\in Q$. 
But then we do not know whether $\{E_a\st a\in Q\}$ is 
order-isomorphic to $Q$, 
because we cannot apply Lemma~\ref{lem:pdomn} 
to prove $E_a\subseteq E_b$ for $a,b\in Q$ with $a<b$ and $\rk{a}=\rk{b}$. 
\end{remark}

%\section{Question}

%\begin{q}
%Does the statement of Hechler's theorem hold 
%for the null ideal of the real line with respect to set-inclusion? 
%\end{q}

\section*{Acknowledgement}

The second author 
would like to thank Prof.\ J\"org Brendle and Mr.\ Teruyuki Yorioka 
for helpful discussion concerning this work.

%\bibliographystyle{plain}
%\bibliography{kada}

\begin{flushleft}
Tomek Bartoszy\'nski \\
Department of Mathematics, Boise State University \\
1910 University Drive, Boise, Idaho 83725 USA\\
E-mail: \texttt{tomek@math.boisestate.edu}
\end{flushleft}

\begin{flushleft}
Masaru Kada \\
Department of Computer Sciences, Kitami Institute of Technology \\
165 Koen-cho, Kitami, Hokkaido 090-8507 JAPAN \\
E-mail: \texttt{kada@math.cs.kitami-it.ac.jp}
\end{flushleft}

\end{document}